\documentclass[final]{siamltex}

\usepackage[parfill]{parskip}    
\usepackage{amsmath,amssymb,latexsym,enumerate,mathrsfs}
\usepackage{hyperref}
\usepackage{array}
\usepackage{times}
\usepackage{subfig}
\usepackage[]{mdframed}
\usepackage{epstopdf}
\usepackage{subfig}
\usepackage{graphicx}
\usepackage{multicol}
\usepackage{moreverb}
\usepackage{marvosym}
\usepackage{color,soul}
\definecolor{lightblue}{rgb}{.90,.95,1}
\sethlcolor{yellow}

\DeclareMathOperator*{\argmin}{arg\,min}

\newtheorem{remark}{Remark}

\newtheorem{assumption}{\bf Assumption}

\title{Embedding generalization within the learning dynamics: An approach based-on sample path large deviation theory}

\author{Getachew K. Befekadu \thanks{This work was partly performed, while the author was visiting Osaka University, Japan.}}

\begin{document}
\maketitle

\renewcommand{\thefootnote}{\arabic{footnote}}

\begin{abstract}
In this paper, we consider a typical learning problem of point estimations for modeling of nonlinear functions or dynamical systems in which generalization, i.e., verifying a given learned model or estimated parameters, can be embedded as an integral part of the learning process or dynamics. In particular, we consider an empirical risk minimization based learning problem that exploits gradient methods from continuous-time perspective with small random perturbations, which is guided by the training dataset loss. Here, we provide an asymptotic probability estimate in the small noise limit based-on the Freidlin-Wentzell theory of large deviations, when the sample path of the random process corresponding to the randomly perturbed gradient dynamical system hits a certain target set, i.e., a rare event, when the latter is specified by the testing dataset loss landscape. Interestingly, the proposed framework can be viewed as one way of improving generalization and robustness in learning problems that provides new insights leading to optimal point estimates which is guided by training data loss, while, at the same time, the learning dynamics has an access to the testing dataset loss landscape in some form of future achievable or anticipated target goal. Moreover, as a by-product, we establish a connection with optimal control problem, where the target set, i.e., the rare event, is considered as the desired outcome or achievable target goal for a certain optimal control problem, for which we also provide a verification result reinforcing the rationale behind the proposed framework. Finally, we present a computational algorithm -- a two-step iterative numerical scheme -- that solves the corresponding variational problem, i.e., a large deviation minimizer problem, leading to an optimal point estimates and, as part of this work, we also present some numerical results for a typical case of nonlinear regression problem.

{\bf Significance.} In most learning problems, the tenet of learning is generally guided by the training dataset loss, while the performance of generalization is evaluated routinely on different novel test datasets, which, in turn, makes it difficult to connect meaningfully what properties of a given learned model will validate the generalizability without having an access to the testing dataset loss landscape. Here, we claim that it is possible to embed generalization in a typical learning framework provided that the learning process or dynamics has access to the testing dataset loss landscape in some form of future achievable or anticipated target goal. To be more specific, there are two main integral parts of the proposed framework. (i) First, we provide an asymptotic probability estimate in the small noise limit based-on the Freidlin-Wentzell theory of large deviations, when the sample path of the random process corresponding to the randomly perturbed gradient dynamical system, which is guided by the training dataset loss, is ultimately nudged towards a certain target set, i.e., a rare event, where the latter is specified by the testing dataset loss landscape. This further allows us to establish a large deviation principle to a variational problem, i.e., a large deviation minimizer problem, modeling the most likely sample path of the rare event leading to optimal point estimates. (ii) Then, with this established large deviation principle result, we provide a computational algorithm that solves the corresponding variational problem. As a by-product, we also establish a connection with optimal control problem, where the target set, i.e., the rare event, is considered as the desired outcome for a certain optimal control problem.

\end{abstract}
\vspace{-0.1in}

\section{Introduction} \label{S1}
In most learning problems, the tenet of learning is generally guided by the training dataset loss, while the performance of generalization is evaluated routinely on different novel test datasets, which, in turn, makes it difficult to connect meaningfully what properties of a given learned model or estimated parameters will validate the generalizability without having an access to the testing dataset loss landscape. In this paper, without attempting to give a literature review, we present a typical learning problem of point estimations for modeling of nonlinear functions or dynamical systems in which generalization, i.e., verifying a given learned model or estimated parameters, can be embedded as an integral part of the learning process or dynamics  (e.g., see \cite{r1a} for general discussions on generalization in recent machine learning literature). Here, we specifically consider an empirical risk minimization based learning problem that exploits gradient methods, which is guided by the training dataset loss, from continuous-time perspective with small random perturbations. Moreover, we provide an asymptotic probability estimate in the small noise limit case based-on the Freidlin-Wentzell theory of large deviations, when the sample path of the random process corresponding to the randomly perturbed gradient dynamical system hits a certain target set, i.e., a rare event, when the latter is specified by the testing dataset loss landscape. 

Interestingly, the proposed framework can be viewed as one way of improving generalization and robustness in learning problems that provides new insights leading to optimal point estimates which is guided by training data loss, while, at the same time, the learning dynamics has an access to the testing dataset loss landscape in some form of future achievable or anticipated target goal. To be more specific, there are two main integral parts of the proposed framework. (i) First, we provide an asymptotic probability estimate in the small noise limit case based-on the Freidlin-Wentzell theory of large deviations, when the sample path of the random process corresponding to the randomly perturbed gradient dynamical system, which is guided by the training dataset loss, is nudged to hit a certain target set, i.e., a rare event, where the latter is specified by the testing dataset loss landscape. This further allows us to establish a large deviation principle to a variational problem, i.e., a large deviation minimizer problem, modeling the most likely sample path of the rare event leading to optimal point estimates. (ii) Then, with this established large deviation principle result, we provide a computational algorithm -- based on a two-step iterative numerical scheme -- that solves the corresponding variational problem. Finally, as a by-product, we also establish a connection with optimal control problem, where the target set, i.e., the rare event, is considered as the desired outcome for a certain optimal control problem, for which we provide verification reinforcing the rationale behind the proposed framework.

This paper is organized as follows. In Section~\ref{S2}, we briefly discuss an empirical risk minimization based learning problem that exploits gradient methods from continuous-time perspective with small random perturbations. Here, we also present some preliminary results based on the Freidlin-Wentzel theory of large deviations that provides a new insight for understanding of learning process corresponding to such gradient systems with small random perturbations. In Section~~\ref{S3}, we present our main results, where we characterize the asymptotic behavior for such small noise diffusions corresponding to continuous-time gradient systems with small random perturbations, and we also provide a rare event interpretation for such learning dynamics with respect to both training and testing datasets loss landscapes.This section also provides a rare event algorithm, i.e., a numerical scheme to compute the large deviation minimizer problem, that characterizes the  maximum likelihood pathways for such learning process. In Section~\ref{S4}, we also present some numerical results for a typical case of nonlinear regression problem, and  Section~\ref{S5} contains concluding remarks.

\section{Preliminaries}\label{S2}
Consider a typical empirical risk minimization based learning problem of point estimators, e.g., a classical case of nonlinear regression problem, where the problem statement is: Given the dataset $Z^n = \bigl\{ (x_i, y_i)\bigr\}_{i=1}^n$, then search for the parameters $\theta \in \Gamma \subset \mathbb{R}^p$, such that $h_{\theta}(x) \in \mathcal{H}$ from a given class of hypothesis function space, that describes best the given dataset. That is, in terms of mathematical optimization construct, we have the following optimization problem
\begin{align}
 \min_{\theta \in \Gamma \subset \mathbb{R}^p} ~ J(\theta, Z^n), \label{Eq2.1}
\end{align}
where
\begin{align}
 J(\theta, Z^n) = \frac{1}{n} \sum_{i=1}^n {\ell}\bigl(h_{\theta}(x_i), y_i\bigr) + \lambda \mathcal{R}(\theta), \label{Eq2.2}
\end{align}
$\ell$ is a suitable loss function that quantifies the lack-of-fit between the model and the dataset, $\lambda > 0$ is regularization parameter that controls the amount of regularization, $\mathcal{R}$ is regularization functional, and $\Gamma \subset \mathbb{R}^p$ is a finite dimensional parameter space. 

In general, the choice for the best regularization is concerned with compromising among different conflicting requirements such as {\it small bias} and {\it small variance}, or determining how faithful the modified problem is to the original problem with better generalization and robustness or smoother optimization decision boundaries (e.g., see \cite{r1} and \cite{r2} for related discussions). On the other hand, if we have prior knowledge via additional constraints on some of the output variable(s) in the form of $\mathcal{Y}(f_{\theta}(x)) \le \eta$. Then, the learning problem can be reformulated as follows
\begin{align}
\min_{\theta \in \Theta} ~ \frac{1}{n} \sum_{i=1}^n \ell\bigl(h_{\theta}(x_i), y_i\bigr) ~~  \text{s.t.} ~~ \mathcal{Y}(f_{\theta}(x)) \le \eta. \label{Eq2.3}
\end{align}
Note that the minimum of $J(\theta, Z^n) \colon \mathbb{R}^p \to \mathbb{R}$ can be located from the asymptotic behavior of the solution corresponding to the following gradient dynamical system which is guided by the dataset $Z^n$
\begin{align} 
 \dot{\theta}(t) = - \nabla J(\theta(t), Z^n), \quad \theta(0) = \theta_0,  \quad  \left( \text{with} ~~ \dot{\theta}(t) = \tfrac{d \theta(t)}{dt} \right), \label{Eq2.4}
\end{align} 
i.e., the {\it steady-state solution} $\theta(t)\, \to \, \theta_{ss} \equiv \theta^{\ast}$ as $t \, \to \, \infty$. Unfortunately, such an approach has a series short coming due to the solution may become trapped at local minimum of $J(\theta(t), Z^n)$, rather than the global minimum solution. One way of addressing such difficulty is to consider the solutions of the related stochastic differential equations (SDEs) with small random perturbations, i.e.,
\begin{align}
d \Theta_t^{\epsilon} = - \nabla J(\Theta_t^{\epsilon}, Z^n) dt + \sqrt{\epsilon} \, I_p d W_t, \quad \Theta_0^{\epsilon} = \theta_0, \label{Eq2.5}
\end{align} 
where $\epsilon \ll 1$ is a small positive parameter, $I_p$ is a $p \times p$ identity matrix, and $W_t$ is a $p$-dimensional standard Wiener process, and study the asymptotic behavior of the diffusion processes. Note that the SDE in Equation~\ref{Eq2.5} possesses a unique strong solution for each $\theta_0$ on $[0,T]$, i.e., a unique $\mathcal{F}_t$-adapted stochastic process $\Theta_t^{\epsilon}$ such that 
\begin{align}
\mathbb{P} \left\{ \sup_{0 < t <T} \biggl \vert \Theta_t^{\epsilon}  - \left( \theta_0 - \int_0^t \nabla J(\Theta_s^{\epsilon}, Z^n) ds + \sqrt{\epsilon} \, \int_0^t I_p d W_s \right) \biggl \vert> 0 \right\} = 0. \label{Eq2.6}
\end{align}
This is due to the assumptions that the SDE in Equation~\ref{Eq2.5} satisfies both the Lipschitz and growth conditions for all $t \in [0,T]$ (e.g., see \cite{r3a} for additional discussions).

\begin{assumption}\label{AS1}
Throughout this work, we assume that $J(\theta,Z^n)$ satisfies the following two conditions:
\begin{enumerate} [(i)] 
\item {\it Coercivity or superlinear growth condition:}
\begin{align}
 \lim_{\theta \to \infty} \frac{J(\theta, Z^n)}{\vert \theta \vert} \to +\infty. \label{Eq2.7}
\end{align}
\item {\it Tightness condition:}
\begin{align}
\int_{\big\{\theta \colon J(\theta, Z^n) \ge \kappa \big\}} \exp \big(-\tfrac{1}{\epsilon} J(\theta, Z^n) \big) d \theta \le C_{\kappa} \exp\big(-\kappa/ \epsilon\big), \label{Eq2.8}
\end{align}
where the constant $C_{\kappa}$ depends on $\kappa$, but not on $\epsilon$.
\end{enumerate}
\end{assumption}
Note that the above two conditions are required to confine the diffusion process $\Theta_T^{\epsilon}$ in a certain compact region with high probability. Here, we also remark that on the consistency of the learning dynamics, where the solution for the SDE in Equation~\ref{Eq2.5} is assumed to be continuously depend on the initial condition $\theta_0$ as well as on the problem dataset $Z^n$. The following general result is useful for verifying further the {\it consistency} of the learning process based on gradient dynamical systems methods with small random perturbations (see also the condition in Equation~\ref{Eq2.6}). 

\begin{proposition}
Suppose that the dataset is perturbed, in the sense that, $Z^n \to Z^{n,\delta} \in \mathcal{D}\bigl(Z^n\bigr) = \left \{ Z^{n,\delta} \, \bigl\vert \rho\bigl(Z^n,Z^{n,\delta} \bigr) \le \delta \right\}$ for some generic distance metric $\rho$ that captures some kinds of perturbation or noises in the dataset.\footnote{e.g., $ \rho\bigl(Z^n,Z^{n,\delta} \bigr) = \bigl\Vert  Z^n - Z^{n,\delta} \bigr \Vert_p$, where $p$ is taken from $\{1,\,2,\, \infty\}$.} Furthermore, assume that $\Theta_t^{\epsilon}$ is a strong solution to the SDE in Equation~\ref{Eq2.5}. Then, we have the following statement 
\begin{align}
\operatorname{st-limit}_{\delta \to 0} \biggl\{ \sup_{t \in [0,T]} \biggl \{ \biggl \vert \Theta_t^{\epsilon,\delta} - \Theta_t^{\epsilon}\biggr \vert \biggr\} \biggr\} = 0, \label{Eq2.9}
\end{align} 
where the stochastic process $\Theta_t^{\epsilon,\delta}$ is a strong solution to the following SDE
\begin{align}
d \Theta_t^{\epsilon,\delta} = - \nabla J(\Theta_t^{\epsilon,\delta}, Z^{n,\delta}) dt + \sqrt{\epsilon} \, I_p d W_t^{\delta}, \quad \Theta_0^{\epsilon,\delta} = \theta_0 \label{Eq2.10}
\end{align}  
and $W_t^{\delta}$ is a standard Wiener process which is independent to $W_t$.
\end{proposition}
\begin{remark}
 The interpretation of the above {\it stochastic-limit} is that the probability of the maximum deviation between $\Theta_t^{\epsilon,\delta}$ and $\Theta_t^{\epsilon}$ over any finite interval $[0,T]$ is nonzero, but it goes to zero as $\delta \to 0$.
\end{remark}

\begin{center}
\begin{figure}[ht]
\begin{center}
 \includegraphics[scale=0.7]{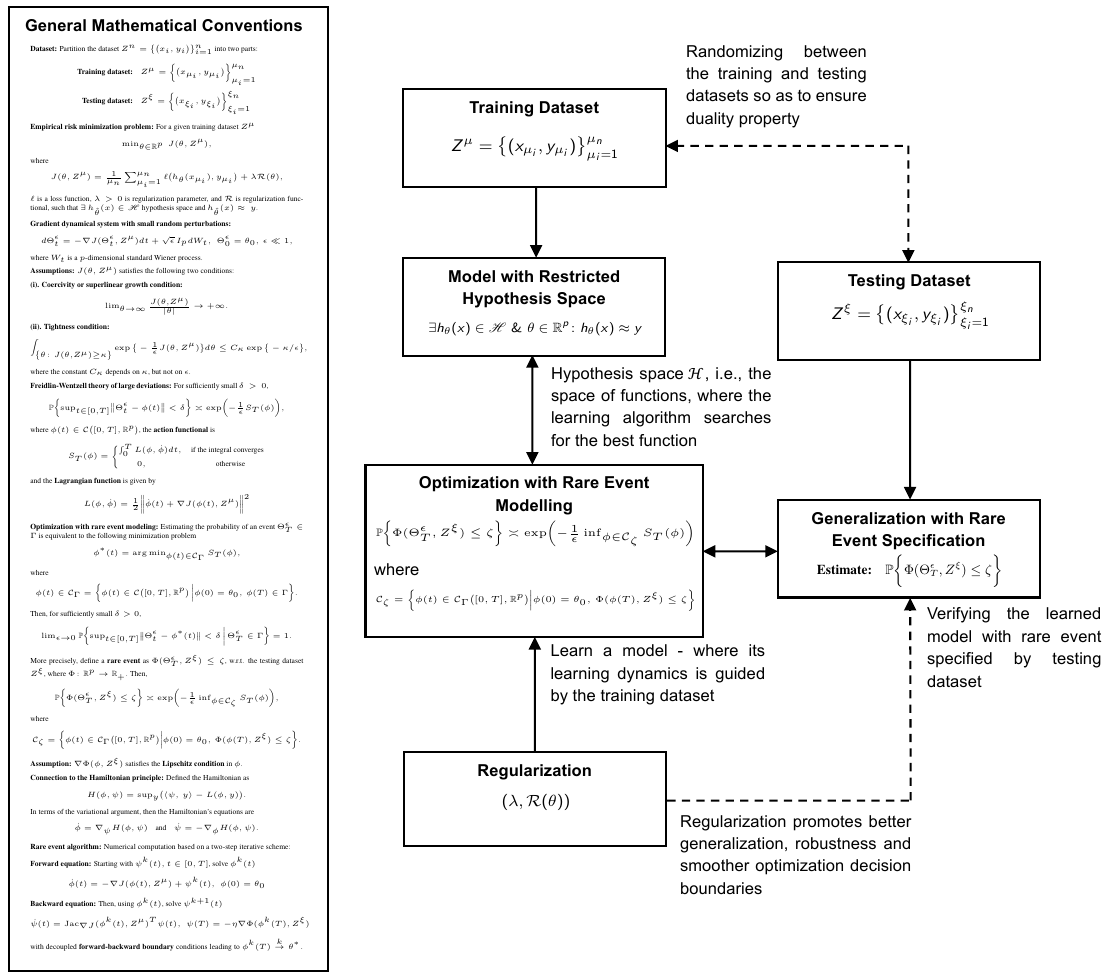}
  \caption{ Flowchart of the framework} \label{FG1}
\end{center}
\end{figure}
\end{center}

In what follows, we partition the dataset $Z^n = \bigl\{ (x_i, y_i)\bigr\}_{i=1}^n$ into two datasets: (i) The first dataset, i.e., the training dataset $Z^\mu = \bigl\{ (x_{\mu_i}, y_{\mu_i})\bigr\}_{\mu_i=1}^{\mu_n}$, will be used for guiding the learning dynamics or model fitting process. (ii) The second dataset, i.e., the testing dataset $Z^\xi = \bigl\{ (x_{\xi_i}, y_{\xi_i})\bigr\}_{\xi_i=1}^{\xi_n}$ which is related to generalization and robustness of the learning problem, will be used for evaluating the performance of the learned model or estimated parameters via specification of rare event corresponding to small noise diffusions of gradient dynamical systems. Figure~\ref{FG1} shows the flowchart of the proposed framework and some of the mathematical conventions that are adopted in this paper. Later in Section~\ref{S3}, we provide the rationale behind the claim of the proposed framework in which it is possible to embed generalization in a typical learning framework provided that the learning process has access to the testing dataset loss landscape in some form of future achievable or anticipated target goal, where such a learning framework will allow us to establish a large deviation principle to a variational problem modeling the most likely sample path of the rare event leading to optimal point estimates. With this established large deviation principle result, we will be able to provide a computational algorithm that solves the corresponding variational problem.

\section{Main results}\label{S3}
In this section, we present our main results, where we present a mathematical formalism, based on Freidlin-Wentzell theory of large deviations and rare event modeling, relevant for understanding and improving generalization in learning problems in which the asymptotic behavior of the learning dynamics provides a new insight leading to optimal point estimates, which is guided by training data loss, while, at the same time, it has an access to the testing dataset loss landscape in some form of future achievable or anticipated target goal. 

\subsection{Freidlin-Wentzell theory of large deviations}\label{S3.1}
Consider the following gradient dynamical system with small random perturbations (cf. the SDE in Equation~\ref{Eq2.5})
\begin{align}
d \Theta_t^{\epsilon} = - \nabla J(\Theta_t^{\epsilon}, Z^\mu) dt + \sqrt{\epsilon} \, I_p d W_t,  \label{Eq3.1}
\end{align}
where $\epsilon \ll 1$ is a small positive parameter, $I_p$ is a $p \times p$ identity matrix, $W_t$ is a $p$-dimensional Wiener process, and $Z^\mu = \bigl\{ (x_{\mu_i}, y_{\mu_i})\bigr\}_{\mu_i=1}^{\mu_n}$ is the training dataset. 

Here, we are mainly interested in situation, where the stochastic process in Equation~\ref{Eq3.1} realizes certain events, e.g., when the trajectory $ \Theta_t^{\epsilon}$ ends in a given set $\Gamma \subset \mathbb{R}^p$, so that $\Theta_T^{\epsilon} \in \Gamma$ at a time $T$. Note that even if these events are impossible in the deterministic system (i.e., when $\epsilon = 0$). They will, in general, occur in the presence of noises (i.e., when $\epsilon > 0$), but they become rarer and rarer in the low-noise limits as $\epsilon \, \to \, 0$. Here,  our analysis mainly relies on the Freidlin-Wentzel theory of large deviations that gives a precise characterization of this decay probability corresponding to learning process or dynamics (e.g., see \cite{r3}, \cite{r4}, \cite{r6}, \cite{r6b} \cite{r7} or \cite{r8} for additional discussions on large deviation principle for such small noise diffusions). Hence, our focus is on an indirect study of learning dynamics based on the gradient methods from the continuous-time perspective with small random perturbations, when the asymptotic behavior of the corresponding diffusion processes and their probabilistic interpretations are exactly known (e.g., see \cite{r3}, \cite{r4} and \cite{r5} for related discussions).

Note that the probability of observing any sample paths close to a given continuous function $\phi(t) \in \mathcal{C} \bigl([0,\,T],\, \mathbb{R}^p \bigl)$ can be estimated as follows:
\begin{align}
\mathbb{P} \biggl \{  \sup\limits_{ t \in [0,\, T]} \biggl \Vert  \Theta_t^{\epsilon} - \phi(t) \biggr \Vert < \delta \biggr\}  \asymp \exp\biggl(-\frac{1}{\epsilon} S_T(\phi)\biggr)  \label{Eq3.2}
\end{align}
for sufficiently small $\delta$, where the notation $\asymp$ denotes {\it log-asymptotic equivalent}, i.e, the ratio of the logarithms of both sides converges to one, and the {\it action functional} $S_T(\phi)$ is given by
\begin{equation}
S_T(\phi) = \left\{ \begin{matrix} 
\int_0^T L(\phi, \dot{\phi})dt, \quad \text {if the integral converges} \\
0, \quad\quad\quad\quad\quad\quad\quad\quad\quad \text {otherwise,} 
\end{matrix}\right.  \label{Eq3.3}
\end{equation}  
where the {\it Lagrangian function} $L(\phi, \dot{\phi})$ is given by
\begin{align}
L(\phi, \dot{\phi}) = \frac{1}{2} \biggl \Vert  \dot{\phi}(t) +  \nabla J(\phi(t), Z^\mu) \biggr \Vert^2,  \label{Eq3.4}
\end{align}
where $\Vert \cdot \Vert^2$ denotes the Riemannian norm of a tangent vector.

In particular, the probability of observing the random event $\Theta_T^{\epsilon} \in \Gamma$ consists of contributions of sample paths close to all possible absolutely continuous function $\phi(t) \in \mathcal{C}_{\Gamma}$, i.e.,
\begin{align}
\phi(t) \in \mathcal{C}_{\Gamma} = \biggl\{ \phi(t) \in \mathcal{C} \bigl([0,\,T],\, \mathbb{R}^p \bigl) \, \bigg\vert \, \phi(0) = \theta_0, \, \phi(T) \in \Gamma \biggr\}.  \label{Eq3.5}
\end{align}
Moreover, in the limit as $\epsilon \to 0$, this is the only contribution which is significantly coming from the trajectory $\phi^{\ast}(t)$ with the {\it smallest action} $S_T(\phi^{\ast})$, i.e.,  
\begin{align}
\phi^{\ast}(t)  = \argmin\limits_{ \phi(t) \in \mathcal{C}} S_T(\phi),  \label{Eq3.6}
\end{align}
which is the maximum likelihood pathway, i.e., the instanton. That is, it constitutes {\it almost surely} all sample paths conditioned on the rare event which is arbitrarily close to $\phi^{\ast}(t)$, while the functional $S_T$ effectively characterizes the difficulty of the passage of $\Theta_t^{\epsilon}$ near $\phi^{\ast}(t)$ in the interval $[0, t]$.

More precisely, for sufficiently small $\delta >0$, we have the following result
\begin{align}
\lim\limits_{\epsilon \to 0} \mathbb{P} \biggl \{ \sup\limits_{ t \in [0,\, T]} \biggl \Vert  \Theta_t^{\epsilon} - \phi^{\ast}(t) \biggr \Vert < \delta \,\biggl\vert \,\Theta_T^{\epsilon} \in \Gamma \biggr\} = 1.  \label{Eq3.7}
\end{align}
In this paper, one of the focus is to provide a computational scheme for the minimizer problem in Equation~\ref{Eq3.5} based on rare events simulations associated with the underlying learning dynamics.

\subsection{Connection with the Hamiltonian principle}\label{S3.2}
The minimization problem in Equation~\ref{Eq3.5}, i.e., finding the maximum likelihood pathway corresponds to the Hamiltonian principle from classical mechanics. That is, the corresponding variational problem can be solved by seeking solutions corresponding to the Euler-Lagrange equation
\begin{align}
\frac{\partial L(\phi,\dot{\phi})}{\partial \phi} - \frac{d}{dt} \frac{\partial L(\phi,\dot{\phi})}{\partial \dot{\phi}} = 0.  \label{Eq3.8}
\end{align}
Next, define the conjugate moment as follows
\begin{align}
\psi = \frac{\partial L(\phi,\dot{\phi})}{\partial \dot{\phi}}.  \label{Eq3.9}
\end{align}
Then, we can define the Hamiltonian as the Fenchel-Legendre transform of the Lagrangian function
\begin{align}
H(\phi,\psi) = \sup\limits_{y} \bigl(\langle \psi,\, y\rangle - L(\phi,y) \bigr)  \label{Eq3.10}
\end{align}
such that the Lagrangian (assuming convexity of $L(\phi, \dot{\phi})$ in $\dot{\phi}$ and noting the duality relation) can be expressed as
\begin{align}
L(\phi,\dot{\phi}) = \sup\limits_{\psi} \bigl(\langle \psi,\, \dot{\phi} \rangle - H(\phi,\psi) \bigr).  \label{Eq3.11}
\end{align}

Consequently, the minimization problem in Equation~\ref{Eq3.5} is equivalent to solving the following Hamiltonian equations of motion, with appropriate boundary conditions,
\begin{align}
\dot{\phi} = \nabla_{\psi} H(\phi,\psi) \quad \text{and} \quad \dot{\psi} = -\nabla_{\phi} H(\phi,\psi),  \label{Eq3.12}
\end{align}
where $H(\phi,\psi) = \bigl \langle - \nabla J(\phi, Z^\mu), \, \psi \bigr \rangle + \tfrac{1}{2}\langle \psi, \psi \rangle$, i.e., 
\begin{align}
\dot{\phi}(t) = - \nabla J(\phi(t), Z^\mu)  + \psi(t) \quad \text{and} \quad  \dot{\psi}(t) = \operatorname{Jac}_{\nabla J}(\phi(t), Z^\mu)^T \psi(t)  \label{Eq3.13}
\end{align}
and $\operatorname{Jac}_{\nabla J}(\phi(t), Z^\mu )$ is the Jacobian of $\nabla J(\phi, Z^\mu)$.

\subsection{Optimization with rare event modeling}\label{S3.3}
Consider the following event that we associate with a rare event
\begin{align}
\Phi\bigl(\Theta_T^{\epsilon}, Z^\xi\bigr) \le \zeta,  \label{Eq3.14}
\end{align}
where $\zeta$ is very small positive number,  $Z^\xi = \bigl\{ (x_{\xi_i}, y_{\xi_i})\bigr\}_{\xi_i=1}^{\xi_n}$ is the testing dataset, and $\Phi\colon \mathbb{R}^p \to \mathbb{R}_{+}$, in the limit as $\epsilon \to 0$, i.e., the probability of observing the event $\Phi\bigl(\Theta_T^{\epsilon}, Z^\xi\bigr) \le \zeta$, subject to $\Theta_0^{\epsilon}=\theta_0$, satisfies the following 
\begin{align}
\mathbb{P} \bigl( \Phi\bigl(\Theta_T^{\epsilon}, Z^\xi\bigr) \le \zeta \bigr) \asymp \exp\biggl ( -\frac{1}{\epsilon} \inf\limits_{\phi \in \mathcal{C}_{\zeta}} S_T(\phi) \biggr),  \label{Eq3.15}
\end{align}
where $\mathcal{C}_{\zeta} = \bigl\{ \phi(t) \in \mathcal{C}_{\Gamma} \bigl([0,\,T],\, \mathbb{R}^p \bigl) \, \big\vert \, \phi(0) = \theta_0, \, \phi(T) \le \zeta \bigr\}$. 

Notice that the condition $\phi(T) \le \zeta$ (or, similarly the random event $\Phi\bigl(\Theta_T^{\epsilon}, Z^\xi\bigr) \le \zeta$) can be considered as some form of future achievable or anticipated target goal, where the testing dataset loss landscape is embedded as an integral part of the learning dynamics or process. Moreover, we assume such a condition as a functional penalty term in the large deviation minimizer problem for which we can impose the following condition on $\Phi(\phi,Z^{\xi})$.

\begin{assumption}\label{AS2}
$\nabla\Phi(\phi,Z^\xi)$ satisfies the Lipschitz condition in $\phi$.
\end{assumption} 

\begin{remark}
Later, in the numerical simulation part, we specifically consider a function $\Phi$ that signifies model validation w.r.t. the testing dataset $Z^\xi$, when the steady-state solution $\Theta_T^{\epsilon}$ tends the global optimum solution $\theta^{\ast}$.
\end{remark}

Then, noting Assumptions~\ref{AS1} and \ref{AS2}, then we define the following rate function
\begin{align}
I(\zeta) = \inf_{\phi \in \mathcal{C}_{\Gamma}} S_T(\phi)  \label{Eq3.16}
\end{align}
and, moreover, the corresponding  Fenchel-Legendre transform is given as follows
\begin{align}
I^{\ast}(\eta) = \inf_{\phi \in \mathcal{C_{\zeta}}} \bigl(S_T(\phi) - \eta \Phi(\zeta) \bigr).  \label{Eq3.17}
\end{align}
Note that, in terms of the variational argument, the Hamiltonian's equations of motion will become (cf. Equation~\ref{Eq3.12})
\begin{align*}
&\text{Forward Equation:} ~~ \dot{\phi}(t) = - \nabla J(\phi(t), Z^\mu)  + \psi(t), ~~ \phi(0) = \theta_0\\ 
&\text{Backward Equation:} ~~ \dot{\psi}(t) = \operatorname{Jac}_{\nabla J} (\phi(t), Z^\mu)^T \psi(t), ~~ \psi(T) = -\lambda \nabla \Phi(\phi(T),Z^\xi),
\end{align*}
 while the {\it forward-backward boundary conditions} are decoupled.

In what follows, we exploit the forward-backward decoupled boundary conditions in the rare event algorithm that cna be implemented via a two-step of numerical iteration scheme for solving the large deviation minimizer problem.\footnote{that makes use of the St\"{o}rmer-Verlet based numerical method leading to proper long-time behavior (e.g., see \cite{r9}).} In particular, the algorithm consists of the following steps for finding the maximum likelihood pathway.
{\rm \small

{\bf ALGORITHM:}
\begin{itemize}
\item[0.] Start with the ${\rm\,k^{th}}$ guess $\phi^k(t) \in \mathcal{C}_{\zeta}$ for the instanton trajectory. Notice that the algorithm works on the space $\mathcal{C}_{\zeta}$ which is more convenient, rather than on $\mathcal{C}_{\Gamma}$.
\item[1.] Solve the equation
\begin{align*}
 \dot{\psi}(t) = \operatorname{Jac}_{\nabla J} (\phi^k(t), Z^\mu)^T  \psi(t), \quad \psi(T) = -\lambda \nabla \Phi(\phi^k(T),Z^\xi)
 \end{align*}
 backward in time to obtain $\psi^{k}(t)$
 \item[2.] Solve the equation
\begin{align*}
 \dot{\phi}(t) = - \nabla J(\phi(t), Z^\mu)  + \psi^k(t), \quad \phi(0) = \theta_0
 \end{align*}
 forward in time to obtain the next guess  $\phi^{k+1}(t)$
 \item[3.] Iterate until convergence, i.e., $\Vert \phi^{k+1}(t) - \phi^{k}(t)\Vert \le {\rm {\rm tol}}$.
\end{itemize}}

\begin{remark}
Note that, in Step 2 above, the functional $\psi^k(t)$ enters as an additional term to the system equation which nudges the dynamics toward the rare event by reweighing the sample paths which is similar to that of important sampling concept. 
\end{remark}

\subsection{Connection with deterministic optimal control problems}\label{S3.4}
In what follows, using ideas from optimal control theory (e.g., see \cite{r9a} and \cite{r4} for similar arguments from the stochastic control point of view), we present results pertaining to an admissible optimal control $\mathscr{U} \ni u \colon [0,T] \to \mathbb{R}^p$, where $\mathscr{U}$ is a class of $\mathbb{R}^p$-valued continuous functions on $[0,T]$, that the gradient dynamical system (which is guided by the training dataset $Z^\mu$)
\begin{align}
\dot{\theta}(t) = - \nabla J(\theta(t), Z^\mu) + u(t),  \quad  \theta(0)=\theta_0 \label{Eq3.18}
\end{align}
will have the desired target goal $\Phi(\theta(T),Z^\xi) = \zeta$ which is specified by the testing dataset loss landscape. Here, we consider the following optimal control problem
 \begin{align}
& \min_{u \in \mathscr{U}}  V(u),  \label{Eq3.19} \\
                                                                                      & \quad \quad \quad \text{s.t. ~~ Equation}~\ref{Eq3.18}, \notag
\end{align}
where $V(u) = \tfrac{1}{2} \int_0^T \vert u(t) \vert^2 dt$. Note that if we further introduce the following two Lagrange multiplies $\psi \in \mathcal{C}([0,T], \mathbb{R}^p)$ and $\eta \in \mathbb{R}$, then the above optimal control problem in Equation~\ref{Eq3.19} is equivalent to minimizing of the following modified cost functional, i.e.,
\begin{align}
\min ~ \to ~ \bar{V}(u) &= \frac{1}{2} \int_0^T \left\vert u(t) \right\vert^2 dt + \eta\left(\Phi(\theta(T),Z^\xi) - \zeta\right) \notag\\
                                             & + \int_0^T \left \langle \psi(t), \,  \dot{\theta}(t) + \nabla J(\theta(t), Z^\mu) - u(t)\right\rangle dt  + \psi^T(0) \left( \theta(0) - \theta_0 \right). \label{Eq3.20}
\end{align}
From calculus of variations (e.g., see \cite{r9b}), the total variation for the cost functional $\hat{V}(u)$ is given by
\begin{align}
\delta \hat{V}(u) &=  \left \langle u(t) - \psi(t),\, \delta u(t) \right \rangle  + \left \langle \dot{\theta}(t) + \nabla J(\theta(t), Z^\mu) - u(t),\, \delta \psi(t) \right\rangle + \left \langle \theta(0) - \theta_0,\, \delta \psi(0) \right\rangle \notag \\
&+ \eta \left \langle \nabla \Phi(\theta(T),Z^\xi),\, \delta \theta(T) \right\rangle  + \left \langle - \dot{\psi}(t) + \operatorname{Jac}_{\nabla J} (\theta(t), Z^\mu)^T \psi(t),\, \delta \theta(t) \right\rangle  \notag \\
& + \left \langle \psi(T),\, \delta \theta(T) \right\rangle  + \left \langle \delta \psi(0),\, \theta(0) \right \rangle. \label{Eq3.21}
\end{align}
Hence, the necessary condition $\delta \hat{V}(u) = 0$ for the extremum implies the following sufficient optimality conditions
\begin{align}
 \dot{\theta}(t) &= - \nabla J(\theta(t), Z^\mu)  + u(t), \quad \theta(0) = \theta_0 \label{Eq3.22}
 \end{align}
 and
 \begin{align}
\dot{\psi}(t) &= \operatorname{Jac}_{\nabla J} (\theta(t), Z^\mu)^T \psi(t), \quad \psi(T) = -\lambda \nabla \Phi(\theta(T),Z^\xi) \label{Eq3.23}
\end{align}
where $\operatorname{Jac}_{\nabla J}(\theta(t), Z^\mu )$ is the Jacobian of $\nabla J(\theta, Z^\mu)$. Here, we further identify the above conditions, i.e., Equations~\ref{Eq3.22} and \ref{Eq3.23}, are equivalent to the decoupled forward-backward Hamiltonian's equations of motion discussed in Subsection~\ref{S3.2} (cf. Equation~\ref{Eq3.13}, with $\theta(t) \to \phi(t)$ and $u(t) \to \psi(t)$). Moreover, the gradient of the modified cost functional $\hat{V}(u)$ with respect to the optimal control $u(t)$ satisfies the following condition
\begin{align}
\frac{\delta \hat{V}(u)}{\delta u} = u(t) - \psi(t) = 0, \quad  \forall t \in [0, T],  \label{Eq3.24}
\end{align}
which is an optimal noise realization that ultimately forces the gradient dynamical system towards the desired target goal $\Phi(\theta(T),Z^\xi) = \zeta$. Here, we remark that the connection between the optimal control problem in Equation~\ref{Eq3.19} and the proposed framework, where the sample path of the randomly perturbed gradient systems, which corresponds to learning dynamics, hits the target set, i.e., the rare event, provides interesting reinforcement as a by-product for the rationale behind the proposed framework.

\section{Numerical simulation results}\label{S4}
\begin{center}
\begin{figure}[h]
\begin{center}
 \includegraphics[scale=0.5]{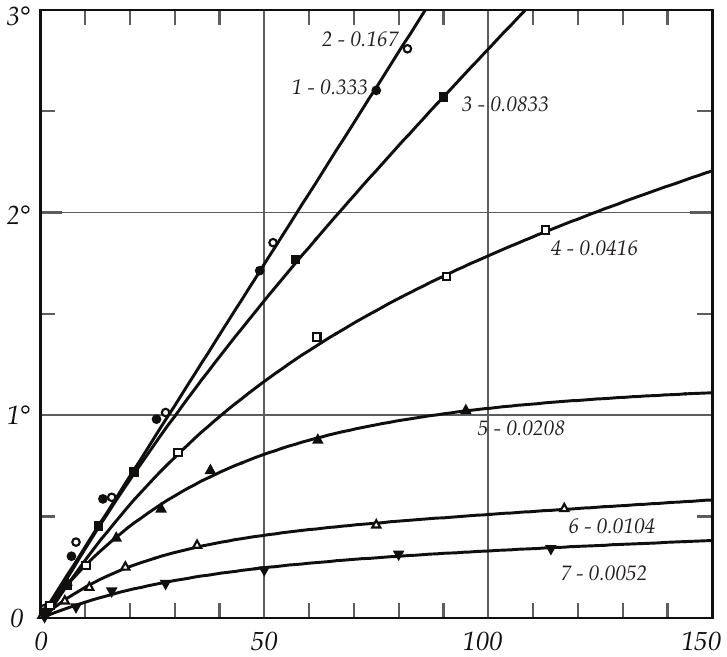}
  \caption{ Abscissa: Time in minutes. Ordinate: Decrease in rotation in degrees.}\label{FG2}
\end{center}
\end{figure}
 \end{center}

Data taken from the Michaelis-Menten's classical paper (see \cite{r18}) that deals with the rate of enzyme-catalyzed reaction. We assume a parametrized functional model, i.e. a nonlinear regression function $v = h_{\theta}(w)$, that relates the initial velocity $v$ as a function of sucrose concentration $w$ from seven different experiments (see also Table~\ref{Tb1}) as follows
\begin{align*}
h_{\theta} (w) &= \frac{\theta_0 \, w}{\theta_1 + w} \quad \text{or} \quad h_{\theta}(w) = \frac{\theta_0 \, w}{\gamma \theta_1 + w}, \quad \text{with $\gamma > 0$ scaling parameter},
\end{align*}
where $\theta_0$ and $\theta_1$ are parameters to be estimated. Note that each curve from an experiment with the given starting concentration of sucrose is shown in Fig.~\ref{FG2}. Here, our interest is to characterize precisely the learning dynamics based on the asymptotic behavior that corresponds to continuous-time gradient dynamical systems with small random perturbations.
\begin{table}[h]
\begin{center}

\caption{Results from the experiments.} \label{Tb1}
  \begin{tabular}{c|c|c}\hline \hline
   &$w$  & $v$ \\
    {\rm Experiments} &{\rm Initial Concentration} & {\rm Initial Velocity} \\ 
    & {\rm  of Sucrose} & \\\hline
    1 &  0.3330 & 3.6360  \\
    2 &  0.1670 & 3.6360 \\
    3 & 0.0833 & 3.2360 \\
    4 &  0.0416 & 2.6660 \\
    5 &  0.0208 & 2.1140 \\
    6 &  0.0104 & 1.4660 \\
    7 &  0.0052 & 0.8661\\
    \hline \hline
  \end{tabular}
  \end{center}
 \end{table}

 In what follows, we presented some numerical results for the above nonlinear functional regression problem based on:
\begin{enumerate} [(i)]
\item {\it Direct ensemble simulation:} Here, we considered the gradient dynamical system $\dot{\theta}(t) = - \nabla J(\theta(t), Z^n)$ with/without small random perturbations, which is guided by the whole dataset $Z^n$, i.e., a few number of experiments $n = 7$. Here, we specifically used a simple Euler--Maruyama time discretization method to solve the system of SDEs (e.g., see \cite{r19} or \cite{r20} for additional discussions on numerical solution of SDEs). Fig.~\ref{FG3} shows the histogram plots with/without noise level $\epsilon = 0.0001$, and bootstrapping the dataset with replacement (i.e., bootstraps of $N = 350$) (e.g., see \cite{r10} for recent discussions on learning via dynamical systems). Moreover, Table~\ref{Tb2} contains the corresponding sample means and variances for direct ensemble simulation.
\begin{figure}[h]
\begin{center}
\includegraphics[scale=0.11]{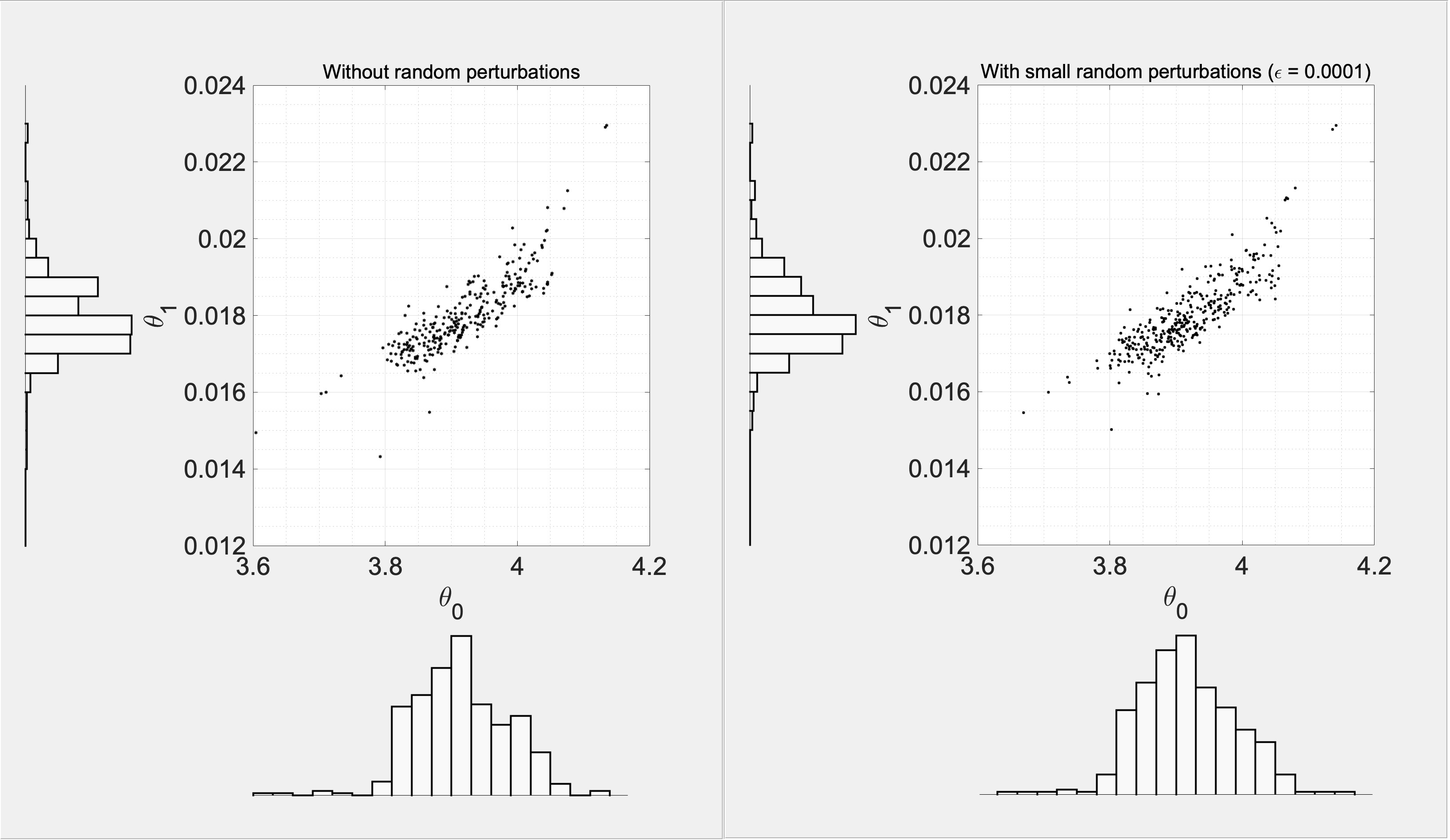} 
\caption{Histograms for the steady-state solutions ($N=350$ bootstraps).} \label{FG3}
\end{center}
\end{figure}
{\normalsize
\begin{table}[h]
\begin{center}
\caption{Sample means and variances for direct ensemble simulation ($N=350$ bootstraps).} \label{Tb2}
\begin{tabular}{l|c|c} \hline \hline 
  Parameters &$\theta_0$  & $\theta_1$ \\ \hline  
   Sample means:  $\bar{\Theta}_N = \tfrac{1}{N}\sum_{i=1}^N \Theta_T^{\epsilon, i}$ &  ${3.9161}$  & ${ 0.0179 }$\\  \hline
   Sample variances: $S_{N-1}^2 =\tfrac{1}{N-1}\sum_{i=1}^N \bigl(\Theta_T^{\epsilon, i} - \bar{\Theta}_N \bigr)^2$ & ${5.3 \times 10^{-3}}$ & ${ 1.9 \times 10^{-6}}$\\ \hline \hline 
  \end{tabular}
  \end{center}
  \end{table}}
\item {\it Rare event characterization via large deviations:} Here, we considered the following rare event that characterizes how near the solution to the optimal parameters value
\begin{align*}
\Phi\bigl(\Theta_T^{\epsilon},Z^\xi\bigr) \le \zeta = 1.5 \times 10^{-5},
\end{align*} 
where the random event $\Phi\bigl(\Theta_T^{\epsilon},Z^\xi\bigr)= (1/\xi_n) \sum\nolimits_{\xi_i=1}^{\xi_n} \bigl(h_{\Theta_T^{\epsilon}}(w_i) - v_i\bigr)^2$ which signifies model validation at a level of order $1.5 \times 10^{-5}$ with limited dataset, when the steady-state solution $\Theta_T^{\epsilon} \to \theta^{\ast} = (3.9109,0.0179)$ that relates the rate for the enzyme-catalyzed reaction. Note that the gradient dynamical system $\dot{\theta}(t) = - \nabla J(\theta(t), Z^\mu)$ with small random perturbations, is required to be guided by the training dataset loss $J(\theta, Z^\mu) = (1/\mu_n) \sum_{\mu_i=1}^{\mu_n} {\ell}\bigl(h_{\theta}(x_i), y_i\bigr)$. Note that the probability of observing the random event $\Phi\bigl(\Theta_T^{\epsilon},Z^\xi\bigr) \le 1.5 \times 10^{-5} $, starting from $\Theta_0^{\epsilon}=\theta_0$, satisfies $\mathbb{P} \bigl(\Phi\bigl(\Theta_T^{\epsilon},Z^\xi\bigr) \le \zeta \bigr) \asymp \exp\bigl (-\tfrac{1}{\epsilon} \inf\nolimits_{\phi \in \mathcal{C}_{\zeta} } S_T(\phi)\bigr)$. Fig.~\ref{FG4} shows additional results based on the rare event and direct ensemble simulations.
\begin{figure}[h]
\begin{center}
   \subfloat[The maximum likelihood pathway, i.e., the instanton trajectory.]{\includegraphics[scale = 0.095]{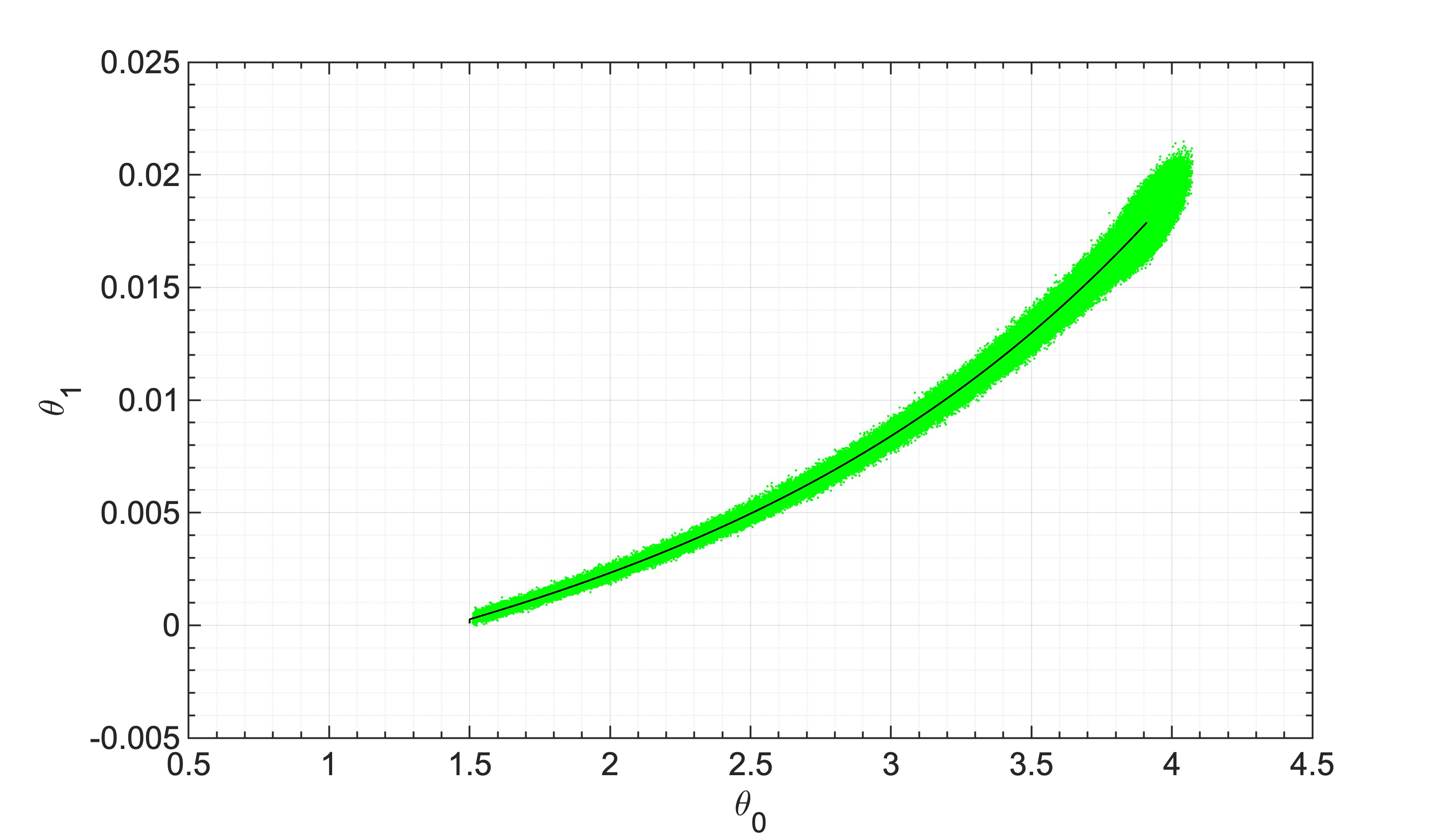}}\\
  \subfloat[Plots for the original dataset, the learned model $h_{\theta^{\ast}}(w)$ and the rare event simulation.]{\includegraphics[scale=0.095]{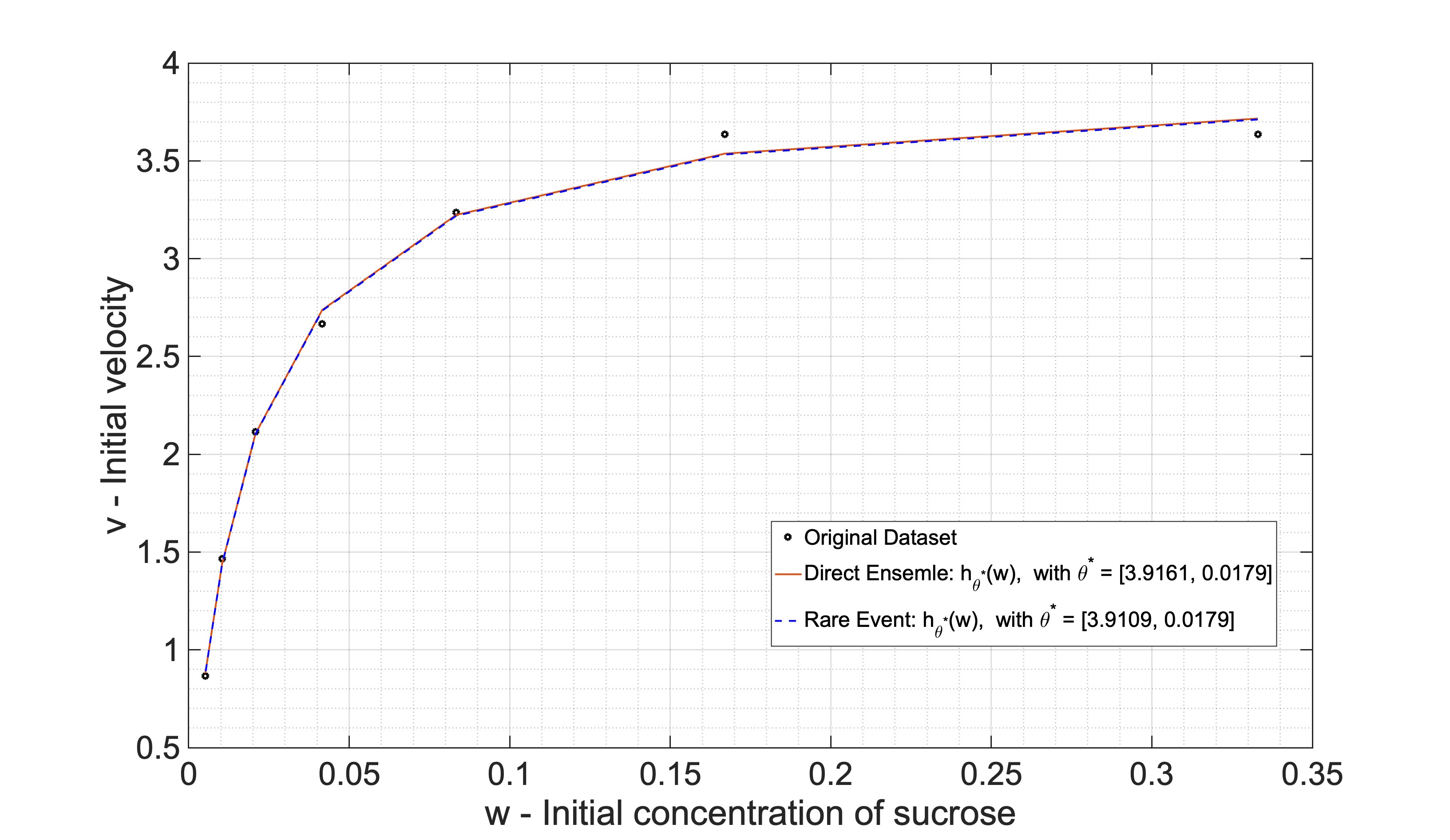}}
    \caption{Results based on the rare event and direct ensemble simulations.} \label{FG4}
    \end{center}
\end{figure}  
\end{enumerate} 

\section{Concluding Remarks}\label{S5}
In this paper, we considered a typical learning problem of point estimations based on empirical risk minimizations for modeling of nonlinear functions in which generalization, i.e., verifying a given learned model or estimated parameters, is embedded as an integral part of the learning process or dynamics. In particular, we presented a mathematical formalism, based on Freidlin-Wentzell theory of large deviations and optimization with rare event modeling, relevant for understanding and improving generalization in learning problems, where the asymptotic behavior of the learning dynamics provides new insights leading to optimal point estimates, which is guided by training data loss, while, at the same time, the learning process has an access to the testing dataset loss landscape in some form of future achievable target goal. Moreover, as a by-product, we also established a connection with optimal control problem, where the target set, i.e., the rare event, is considered as the desired outcome for a certain minimum optimal control problem, for which we provided a verification result reinforcing the rationale behind the proposed framework. Finally, we presented a computational algorithm for solving the corresponding variational problem leading to an optimal point estimates and, as part of this work, we also presented some numerical results for a typical nonlinear regression problem.

\end{document}